\documentclass[12pt]{article} 
\usepackage{amstex, amsthm}

\newcommand{\doublespace}{
   \renewcommand{\baselinestretch}{1.2}
   \large\normalsize}
\renewcommand{\theequation}{\thesection.\arabic{equation}}

\begin{document}

\begin{center}
{\Large {\bf Automorphism Groups and 
Twisted Modules for Lattice Vertex Operator Algebras }} \\
\vspace{0.5cm}
Chongying Dong\footnote{Supported by NSF grant 
DMS-9700923 and a research grant from the Committee on Research, UC Santa Cruz.}\\
Department of Mathematics, University of California\\
Santa Cruz, CA 95064\\
\vspace{0.5cm}
Kiyokazu Nagatomo\footnote{Partly supported by Grant-in-Aid for Scientific Research,
the Ministry of Education, Science and Culture.}
\\
Department of Mathematics, Graduate School of Science\\
Osaka University, Toyonaka, Osaka 560-0048, Japan
\end{center}

\hspace{1.5 cm}
\begin{abstract} We give a complete description of 
the full automorphism group of a lattice vertex operator algebra,
determine the twisted 
Zhu's  algebra for the automorphism $\theta$ 
lifted  from the $-1$ isometry of the lattice
and classify 
the irreducible $\theta$-twisted modules through the twisted Zhu algebra.
\end{abstract}

\renewcommand{\theequation}{\thesection.\arabic{equation}}

\newcommand{\C}{\mathbb{C}}
\newcommand{\Z}{\mathbb{Z}}
\newcommand{\Q}{\mathbb{Q}}
\newcommand{\R}{\mathbb{R}}
\newcommand{\N}{\mathbb{N}}
\newcommand{\1}{\mathbf{1}}
\newcommand{\im}{\operatorname{Im}}
\newcommand{\vac}{|0\,\rangle}
\newcommand{\wt}{\operatorname{wt}\,}
\newcommand{\unit}{\mathbf{1}}
\newcommand{\Hom}{\operatorname{Hom}\,}
\newcommand{\End}{\operatorname{End}}
\newcommand{\LHS}{\operatorname{LHS}}
\newcommand{\RHS}{\operatorname{RHS}}
\newcommand{\Ker}{\operatorname{Ker}}
\newcommand{\frakg}{\mathfrak{g}}
\newcommand{\hatfrakg}{\hat{\mathfrak{g}}}
\newcommand{\eulerz}{z\partial_z}
\newcommand{\eulerw}{w\partial_w}
\newcommand{\id}{\operatorname{id}}
\newcommand{\Vir}{\mathcal{V}ir}
\newcommand{\Id}{\operatorname{Id}}
\newcommand{\ad}{\operatorname{ad\,}}
\newcommand{\Aut}{\operatorname{Aut}\,}
\newcommand{\Outn}{\operatorname{Out}}
\newcommand{\Out}{\operatorname{Out}\,}
\newcommand{\Autn}{\operatorname{Aut}}
\newcommand{\Inn}{\operatorname{Inn}}
\newcommand{\Iso}{\operatorname{Iso}\,}
\newcommand{\NO}{\,{\raise0.25em\hbox{$\mathop{\hphantom{\cdot}}%
\limits^{_{\circ}}_{^{\circ}}$}}\,}

\def \Res{{\rm Res}}
\def \<{\langle}
\def \>{\rangle}
\def \a{\alpha}
\def \b{\beta}
\def \e{\epsilon}
\def \s{\sigma}
\def \o{\omega}

\doublespace

\theoremstyle{plain}
 \newtheorem{theorem}{Theorem}[section]
        \newtheorem{corollary}[theorem]{Corollary}
        \newtheorem{lemma}[theorem]{Lemma}
        \newtheorem{proposition}[theorem]{Proposition}
        \newtheorem{conjecture}[theorem]{Conjecture}
 \newtheorem{definition}[theorem]{Definition}
 \newtheorem{remark}[theorem]{Remark}
 \newtheorem{note}[theorem]{Note}
 \newtheorem{example}[theorem]{Example}
 \newtheorem{exercise}[theorem]{Exercise}

\section{Introduction}
\setcounter{equation}{0}

Vertex (operator) algebra $V_L$ associated with an arbitrary even lattice $L$
is one of the important motivations for introducing the notion of
vertex (operator) algebra and also 
provides basic examples of such algebras
([B], [FLM]). The study of $V_L$ when $L$ is positive has been fruitful.
If $L$ is positive definite $V_L$ is regular [DLM1] in the
sense that any weak module is a direct  sum of ordinary irreducible
modules which have been classified (see [FLM] and [D1]). 
A characterization of $V_L$ was obtained in [LX]. The automorphism
groups of all known vertex operator subalgebras of $V_L$ when the
rank of $L$ is one are determined in [DG] and [DGR]. 

The vertex operator algebra $V_L$ has a canonical 
automorphism $\theta$ of order 2 lifted from the $-1$ isometry
of $L.$ The aim of this paper is to determine the full automorphism group 
$\Aut(V_L)$ of 
$V_L,$ compute the $\theta$-twisted Zhu's algebra [DLM3] and
classify the irreducible $\theta$-twisted modules for $V_L$
by using the twisted Zhu's algebra when $L$ is positive definite.

The group structure of $\Aut(V_L)$ determined in this paper is important
in understanding automorphism groups of vertex operator algebras
in general. Let $V$ be a vertex operator algebra without vectors of
negative weights and weight zero subspace being one dimensional.
Then the weight one subspace $V_1$ is a Lie algebra under the bracket
$[u,v]=u_0v$ where $u_0$ is the component operator of 
vertex operator $Y(u,z)=\sum_{n\in \Z}u_nz^{-n-1}$ associated to
$u.$ The exponentials $e^{u_0}$ for $u\in V_1$ generate a continuous,
normal subgroup $N$ of $\Aut(V).$ So the task is to determine the
quotient group $\Aut(V)/N.$ In this paper we show that 
$\Aut(V_L)/N$ is isomorphic to a quotient group of the finite group 
$O(L)$ which is the subgroup of $\Aut(L)$ of isometries. This may suggest that
$\Aut(V)/N$ is a finite group in general. There are two extremes.
Namely, $V_1$ generates $V$ in the sense of [FHL] or $V_1=0.$ For example,
the vertex operator algebras associated to affine Kac-Moody
algebras (cf. [DL1], [FZ] and [Li]) 
are generated by weight one subspaces and $\Aut(V)$
is isomorphic to the automorphism group of the Lie algebra $V_1;$
 the moonshine vertex operator
algebra $V^{\natural}$ [FLM] and the vertex operator algebras 
associated to the Virasoro algebra (cf. [FZ] and [Li]) 
have zero weight one subspaces
and the automorphism groups are the monster simple group
and the identity. 

The automorphisms studied in this paper preserve the Virasoro
element. But it is also an important problem to 
investigate the automorphisms which need not
preserve the Virasoro element. Such problem was first considered in
[MN] in the case of free bosonic vertex algebra, where  
the the group of all such automorphisms of free
bosonic vertex algebra was obtained. It is still an open problem
to determine all such automorphisms for $V_L.$ 

In [DLM3], an associative algebra $A_g(V)$ was introduced
for a vertex operator algebra $V$ and an automorphism $g$ of
finite order so that there is a one to one correspondence
between the equivalence classes of irreducible admissible
$g$-twisted $V$-modules and the equivalence classes of irreducible
$A_g(V)$-modules. In the case $g=1$ this amounts Zhu's
original algebra $A(V)$ and related results. The Zhu's
algebras $A(V)$ have been computed for the vertex operator
algebras associated to affine Lie algebras in [FZ], 
the vertex operator algebras
associated to the Virasoro algebras in [W] and for $V_L$ in [DLM2].
Our result on $A_\theta(V_L)$ in this paper is the first example
of explicit calculation of twisted Zhu's algebra. It turns out
the computation of $A_\theta(V_L)$ is much simpler than
that of $A(V_L)$ and the structure of $A_\theta(V_L)$ is also simpler.

Let $\hat L$ be the canonical central extension
of $L$ by a cyclic group $\<\pm 1\>$ of order 2. 
In [FLM] and [Le], a subgroup $K$ of $\hat L,$ which is isomorphic
to $2L,$  was introduced and all irreducible modules for $\hat L/K$
on which $(-1)K$ acts as $-1$ were classified and denoted by $T_\chi.$ 
It was shown in [FLM] (also see [DL2]) each $T_\chi$ gives rise an irreducible
$\theta$-twisted module. We show in this paper that $A_\theta(V_L)$
is isomorphic to a quotient of the group algebra $\C[\hat L/K]$ modulo
the ideal generated by $(-1)K+1.$ As a result we 
know all the irreducible modules of $A_\theta(V_L)$ are precisely
these $T_{\chi}.$ Using the theory of twisted Zhu's algebra
developed in [DLM3] we give a new proof of a result in [D2] on 
classification of irreducible $\theta$-twisted $V_L$-module. Although
the proof in [D2] can be easily generalized to prove that $V_L$
is $\theta$-rational but this has not appeared in any publication.
We also prove a general result which says that for an arbitrary 
vertex operator algebra $V$ and an automorphism $g$ of finite order
such that $A_g(V)$ is finite-dimensional, semisimple and the weights of 
top levels of the all irreducible $g$-twisted $V$-modules are the same
then $V$ is $g$-rational. As a corollary we show that 
$V_L$ is $\theta$-rational. 

The most of this work was done when the second author
was visiting UC Santa Cruz in academic year 1997-1998. The second
author wishes to thank the Department of Mathematics here for
their hospitality during his stay. 
The authors are grateful to Geoffrey Mason and Michael Tuite for helpful
comments related to this work. 

\section{Automorphism groups}

In this section we determine the automorphism group of the
vertex operator algebras $V_L$ associated to even positive definite
lattices $L.$ In Subsection 2.1 we review the canonical central extension
$\hat L$ of $L$ by the cyclic group of order 2 and the
structure of the automorphism group of $\hat L$ following [FLM].
We recall the explicit construction of the vertex operator algebra $V_L$
from [FLM] in Subsection 2.2. Subsection 2.3 is on how to lift a certain
derivations of a vertex operator algebra to automorphisms.  Subsection
2.4 is the core of Section 2, where the automorphism group of
$V_L$ is obtained by using the automorphism group of $\hat L$ and
the automorphisms which are the exponentials of the weight zero components
of the weight one vectors. 

\setcounter{equation}{0}
\subsection{Integral lattices and their central extensions}

Let $L$ be any positive-definite even lattice and 
$\hat{L}$ be the canonical central extension of $L$ by the cyclic
group $\< \pm 1\>$:
\begin{equation}\label{2.1}
1\;\rightarrow \< \pm 1\>\;\rightarrow \hat{L}\;\bar{\rightarrow}
L\;\rightarrow 1 
\end{equation}
with the commutator map $c(\alpha,\beta)=(-1)^{\< \alpha,\beta\>}$ for
$\alpha,\beta \in L$. Let $e:  L\to \hat L$ be a section such that
$e_0=1$ and $\epsilon: L\times L\to \<\pm 1\>$ be
the corresponding 2-cocycle. Then $\epsilon(\a,\b)\epsilon(\b,\a)=(-1)^{\<\a,\b\>},$ 
$$\e(\a,\b)\e(\a+\b,\gamma)=\e(\b,\gamma)\e(\a,\b+\gamma)$$
and $e_{\a}e_{\b}=
\e(\a,\b)e_{\a+\b}$ for $\a,\b,\gamma\in L.$
We will use $\Aut(G)$ to denote the automorphism group of any group $G.$ 
Following [FLM] we define 
$$\Aut(c)=\{\sigma\in \Aut(L)|c(\sigma\a,\sigma\b)=c(\a,\b)\ for \ \a,\b\in L\}.$$
Then we have the lifting property
(cf. Proposition 5.4.1 of \cite{FLM})
\begin{equation}\label{2c}
1\rightarrow \Hom (L,\Z/2\Z) \rightarrow
\Aut(\widehat{L})\overset{-}{\rightarrow}
\Aut(c)\rightarrow 1\quad \text{(exact)}.
\end{equation} 
One can easily verify that
\[
\Hom (L,\Z/2\Z)\cong (\Z/2\Z)^{\,\operatorname{rank}L}.
\]

\subsection{Lattice vertex operator algebras}

Let ${\mathfrak h}={\C}\otimes_{{\Z}}L$ and
extend the ${\Z}$-form on $L$ to ${\mathfrak h}$ by $\C$-linearity.
 The corresponding 
affine Lie algebra is $\hat{{\mathfrak h}}={\mathfrak h}\otimes \C[t,t^{-1}]
\oplus \C c.$
We also use the notation $h(n)=h\otimes t^n$ for $h\in {\mathfrak h}, n\in
{\Z}$. Set 
\begin{eqnarray*}
\hat{{\mathfrak h}}^{-}={\mathfrak h}\otimes t^{-1}\C[t^{-1}].
\end{eqnarray*} 
Then $\hat{{\mathfrak h}}^{-}$ is an abelian subalgebra of
$\hat{{\mathfrak h}}$.  
Let $U(\hat{{\mathfrak h}}^{-})=S(\hat{{\mathfrak h}}^{-})$ be
the universal enveloping algebra of $\hat{{\mathfrak h}}^{-}$. Consider the
induced $\hat{{\mathfrak h}}$-module
\begin{eqnarray*}
M(1)=U(\hat{{\mathfrak h}})\otimes _{U(\C[t]\otimes {\mathfrak h}\oplus \C c)}
{\C}\simeq S(\hat{{\mathfrak h}}^{-})\;\;\mbox{(linearly)},\end{eqnarray*}
where $\C[t]\otimes {\mathfrak h}$ acts trivially on $\C$ and $c$ acts
on $\C$ as multiplication by 1.

Recall from (\ref{2.1}) that $\hat{L}$ is the canonical central 
extension of $L$ by the cyclic
group $\< \pm 1\>.$ Form the induced $\hat{L}$-module
\begin{eqnarray*}
{\C}\{L\}=\C[\hat{L}]\otimes _{\< \pm 1\>}\C\simeq
\C[L]\;\;\mbox{(linearly)},\end{eqnarray*}
where $\C[\cdot]$ denotes the group algebra and $-1$ acts on ${\bf
C}$ as multiplication by $-1$. For $a\in \hat{L}$, write $\iota (a)$ for
$a\otimes 1$ in $\C\{L\}$. Then the action of $\hat{L}$ on ${\C}
\{L\}$ is given by: $a\cdot \iota (b)=\iota (ab)$ and $(-1)\cdot \iota
(b)=-\iota (b)$ for $a,b\in \hat{L}$.

Furthermore we define an action of $ {\mathfrak h}$ on $\C\{L\}$ by:
$h\cdot \iota (a)=\< h,\bar{a}\> \iota (a)$ for $h\in {\mathfrak  h},a\in
\hat{L}$. Define $z^{h}\cdot \iota (a)=z^{\< h,\bar{a}\> }\iota (a)$.

The untwisted space associated with $L$ is defined to be
$$V_{L}=M(1)\otimes_\C\C\{L\}$$
and the vertex algebra algebra structure on $V_L$ (see [B] and [FLM])
is determined by 
$$Y(h(-1)\1, z) = h(z)=\sum_{n\in\Z}h(n)z^{-n-1}\quad (h\in\mathfrak{h})$$
$$Y(\iota(a),z) = \exp{(-\sum_{n<0}\frac{\bar a(n)}{n}z^{-n})}
\exp{(-\sum_{n>0}\frac{\bar a(n)}{n}z^{-n})}az^{\bar a}.$$
The Virasoro element is given by
\[
\omega = \frac{1}{2}\sum_{i=1}^dh_i(-1)^2\1
\]
where $\{h_i|i=1,...,d\}$ is an orthogonal basis of $\mathfrak{h}$. 
Expand $L(z) = Y(\omega,z)$
\[
L(z) = \sum_{n\in\Z}L(n) z^{-n-2}.
\]
Then the component operator $L(0)$ defines the weight 
gradation on $V_L:$
\[
V_L = \oplus_{n\in\Z}(V_L)_n,\quad
(V_L)_n = \{v\in V_L\,|\,L(0)v = nv\}.
\]
Note that $\dim\, (V_L)_n<\infty$ for each $n.$ 
For example,
\[
\wt(h(-n)\1) = n,\quad \wt(\iota(a)) =
\frac{1}{2}\langle\bar a,\bar a\rangle.
\]
We should point out that $M(1)$ is a vertex operator subalgebra of
$V_L.$

\subsection{Derivations and automorphisms}

Let $V$ be a vertex operator algebra and $\omega$ the Virasoro
element.  An {\em automorphism} of $V$ is an isomorphism $\sigma:V\longrightarrow V$
of vector spaces which preserves all the products and the Virasoro element:
\[
\sigma(a_nb) = \sigma(a)_n\sigma(b),\quad \sigma(\omega) = \omega.
\]
Note that an automorphism of $V$ preserves each homogeneous subspace
$V_n$ of $V.$
The group of all automorphisms of the vertex operator algebra $V$ is denoted by
$\Aut(V)$. 

A {\em derivation} of $V$ is an endomorphism $D:V\longrightarrow V$ such that
\[
D(a_nb) = D(a)_nb+a_nD(b)\quad \text{for all } a,b\in V, n\in\Z
\]
and $D(\omega) = 0$. In particular, $D$ preserves the gradation of $V.$ 
So the exponential $e^D$ converges on $V$ and is well-defined. 
It is easy to see that $e^D$ is an automorphism.

Let $a\in V$ and $\wt(a)=1$. Then $D= a_0$ is a derivation of $V$ and 
$e^{a_0}$ is an automorphism of $V.$ 
Set $$N=\langle e^{a_0} | a\in V_1\rangle.$$
Since $\sigma e^{a_0}\sigma^{-1} = e^{(\sigma a)_0}$ and 
$\wt(\sigma(a))= 1$ for any $\sigma\in \Aut(V),$ $N$ is a normal subgroup of
$\Aut(V)$. 

\subsection{The automorphism group of $V_L$}

We first describe certain automorphisms of $V_L.$ Set
$$O(L)=\{\sigma\in \Aut(L)|\<\sigma\a,\sigma \b\>=\<\a,\b\>\ for \ \a,\b\in L\}$$
which is the subgroup of $\Aut(L)$ of isometries. Define
$$O(L)=\{\sigma\in \Aut(\hat L)|\bar\sigma\in O(L)\}.$$
Then $O(\hat L)$ is a subgroup of $\Aut(\hat L)$ and $\Hom(L,\Z/2\Z)$
is a subgroup of $O(\hat L)$ (see (\ref{2c})). It follows from (\ref{2c}) that
we have an exact sequence
$$1\rightarrow \Hom (L,\Z/2\Z) \rightarrow
O(\hat{L})\overset{-}{\rightarrow}
O(L)\rightarrow 1.$$

Every $\sigma\in O(\hat L)$ induces an automorphism of $V_L,$ denoting
by $\sigma$ also, in the following way
$$\s(\a_1(-n_1)\cdots \a_k(-n_k)\otimes \iota(a))=(\bar \s\a_1)(-n_1)\cdots (\bar \s\a_k)(-n_k)\otimes \iota(\s a)$$
for $\a_i\in \mathfrak h,$ $n_i>0$ and $a\in \hat L.$ So we can regard
$O(\hat L)$ as a subgroup of $\Aut(V_L).$ Note that such automorphisms
preserve the vertex operator subalgebra $M(1).$

The first main result is the determination of the automorphism group of $V_L.$
\begin{theorem}\label{t2.1}
Let $L$ be a positive definite even lattice. Then
\[
\Aut(V_L) = N \cdot O(\hat L)
\]
and $N$ is a normal subgroup of $\Aut(V_L).$ 
Moreover, the intersection $N\cap O(\hat L)$ contains 
a subgroup $\Hom (L,\Z/2\Z)$ and 
$\Aut(V_L)/N$ is isomorphic to a quotient group of $O(L).$ 
\end{theorem}

To prove this theorem, we need the following elementary lemmas.

\begin{lemma}\label{lemma:2.5}
Let $0\neq v\in V_L$.
Suppose that there exists $\alpha\in \mathfrak{h}\,(\alpha\neq 0)$ such that
$h(0)v = \<h,\alpha\>v$ for all $ h\in\mathfrak{h}$.
Then $\alpha\in L$ and there exists $c_\alpha\in M(1)$ such that
$v =c_\alpha \iota(e_\alpha)$.
\end{lemma}
\begin{proof} Note that $\mathfrak{h}$ acts on $V_L$ semisimply,
$$V_L=\bigoplus_{\beta\in L}M(1)\otimes \iota(e_\beta)$$
and $M(1)\otimes \iota(e_\beta)$ is precisely the subspace with the eigenfunction
$\beta,$ that is 
$$M(1)\otimes \iota(e_\beta)=\{v\in V_L|h(0)v=\left<h,\beta\right>v\}.$$
The lemma follows immediately.
\end{proof}

Set $\mathfrak{h}(-1)=\mathfrak{h}\otimes t^{-1}\1,$ a subspace of $V_L.$ 
Then any automorphism $\sigma$ of $V_L$ preserving $\mathfrak{h}(-1)$
can be considered as a linear automorphism of $\mathfrak{h}.$ 

\begin{lemma}\label{lemma:2.6}
Let $\sigma\in \Aut(V_L)$ and suppose $\sigma(\mathfrak{h}(-1)) =
\mathfrak{h}(-1)$. Then $\sigma (L)\subset L$ and 
for any $\alpha\in L$, 
there exists
$c_\alpha\in\C$ such that
$\sigma(\iota(e_\alpha)) = c_\alpha \iota(e_{\sigma(\alpha)})$.
\end{lemma}
\begin{proof}
From the assumption $\sigma^{-1}(h(-1)\1)\in \mathfrak{h}(-1)$ for any $h\in \mathfrak{h} $. Since the bilinear form on $\mathfrak{h}$ can be realized as
$\<u,v\>=u(1)v$ for $u,v\in \mathfrak{h}$ we see that the form
$\<,\>$ on $\mathfrak{h}$ is $\sigma$-invariant. Here we have identified
$\mathfrak{h}$ with $\mathfrak{h}(-1).$
Then we have from the definition of the action that
$\sigma^{-1}(h)(0)\iota(e_\alpha) = \langle\alpha,\sigma^{-1}(h)\rangle \iota(e_\alpha)$.
Thus
\begin{align*}
h(0)\sigma\iota(e_\alpha)& =
(\sigma\sigma^{-1}(h))(0)\sigma\iota(e_\alpha)\\
& = \sigma(\sigma^{-1}(h)(0)\iota(e_\alpha))\\
&= \sigma(\langle\alpha,\sigma^{-1}(h)\rangle \iota(e_\alpha))\\
& =\langle\alpha,\sigma^{-1}(h)\rangle\sigma\iota(e_\alpha).
\end{align*}
This implies that
$h(0)\sigma\iota(e_\alpha) = \langle h,\sigma(\alpha)\rangle\sigma\iota(e_\alpha)$. 
Clearly $\sigma(\iota(e_\alpha))\neq 0.$ 
Using Lemma \ref{lemma:2.5} we see that
there exists
$c_\alpha\in M(1)$ such that
$\sigma(\iota(e_\alpha)) = c_\alpha \iota(e_{\sigma(\alpha)}).$
In particular, $\sigma (L)\subset L.$ 
Note that $\C\{L\}$ is
the vacuum space for the Heisenberg algebra 
$\oplus_{n\ne 0}t^{n}\otimes \mathfrak{h}\oplus \C c$
in the sense that
$$\C\{L\}=\{v\in V_L|h(n)v=0, h\in  \mathfrak{h}\}.$$
So $\C\{L\}$ is $\sigma$-stable and $c_\alpha\in \C.$ 
\end{proof}

\begin{remark}\label{r1} Since the vectors $\iota(e_{\alpha})$ and
$\iota(e_{\s\alpha})$ in the proof of Lemma \ref{lemma:2.6} 
have the same weight, $\sigma$ induces an isometry of $L.$
\end{remark}

\begin{lemma}\label{lemma:2.7} Let $\sigma\in \Aut(V_L)$ such that
$\sigma|_{\mathfrak{h}(-1)} = \id$. 
Then there exists $h\in\mathfrak{h}$ such that $\sigma =
e^{h(0)}$.
\end{lemma}
\begin{proof}

Since $M(1)$ is generated by $\mathfrak{h}(-1)$ we see immediately that
the restriction of $\sigma$ to $M(1)$ is the identity. Note that
$\sigma$ is also the identity on $L.$ By Lemma \ref{lemma:2.6} for any 
$\alpha\in L$ there exists $c_\alpha\in \C$ so that $\sigma\iota(e_\a)=
c_\a\iota(e_\a).$ Recall the component vertex operators $v_n$
of $Y(v,z)=\sum_{n\in\Z}v_nz^{-n-1}$ for any $v\in V_L.$ Then for 
any $\alpha,\b\in L,$ $\iota(e_\alpha)_{\langle\alpha,\,\beta\rangle-1} \iota(e_\beta)=\epsilon(\alpha,\beta)\iota(e_{\alpha+\beta})$.
We see
\begin{align*}
\sigma\iota(e_{\alpha+\beta}) &= \epsilon(\alpha,\beta)^{-1}
\sigma\iota(e_\alpha)_{\langle\alpha,\,\beta\rangle-1}\sigma\iota(e_\beta)\\
&= \epsilon(\alpha,\beta)^{-1}c_\alpha c_\beta \iota(e_\alpha)_{\langle\alpha,\,\beta\rangle-1} \iota(e_\beta)\\
&= c_\alpha c_\beta \iota(e_{\alpha+\beta}).
\end{align*}
This shows that $c_{\a+\b}=c_\a c_\b.$ 

Let $\{\alpha_i\}_{i=1}^{d}$ be a basis of $L$. Then for $\alpha =
\sum_{i=1}^dk_i\alpha_i\,(k_i\in\Z)$, we see
$c_\alpha = c_1^{k_1} c_2^{k_2}\cdots  c_d^{k_d}$.
Since the inner product is non-degenerate, 
there exists $h\in\mathfrak{h}$ such that
$e^{\langle h,\alpha_i\rangle} = c_{\a_i}$ for all $i=1,2,\dots, d$.
Then  $e^{h(0)}\iota(e_{\alpha_i}) = c_{\alpha_i}\iota(e_{\alpha_i})$
for all $i=1,2,\dots, d$ and 
\[
e^{h(0)}\iota(e_{\alpha}) = c_{\alpha_1}^{k_1} \cdots  c_{\alpha_d}^{k_d} \iota(e_{\alpha})= c_\alpha \iota(e_{\alpha}).
\]
Since both $\sigma$ and $e^{h(0)}$ act trivially 
on $\mathfrak{h}(-1)$ and have the same action on $\C\{L\},$ we conclude that
$\sigma = e^{h(0)}$.
\end{proof}

{\bf Proof of Theorem} \ref{t2.1}. Let us consider the weight one space
$$\mathfrak{g} =V_L(1)= \mathfrak{h}\oplus \sum_{\<\alpha,\alpha\>
= 2}\C \iota(e_\alpha)$$
which is a finite dimensional reductive Lie algebra under the bracket
$[u,v]=u_0v$ for $u,v\in \mathfrak{g}$ with  
a Cartan subalgebra $\mathfrak{h}=\mathfrak{h}(-1).$ 

Let $\sigma$ be an automorphism of $V_L$. Then $\sigma$ induces an automorphism
of the  Lie algebra $\mathfrak{g}.$ Note that the inner automorphism
group of $\mathfrak g$ is the restriction of $N$ to $\mathfrak g.$ 
The conjugacy theorem  for Cartan subalgebras of finite dimensional 
reductive Lie algebras implies that there exists 
$\phi\in N $ such that
$\phi^{-1}\sigma(\mathfrak{h}) =\mathfrak{h}$. 
Then by Lemma \ref{lemma:2.6} and Remark \ref{r1}, $\mu=\phi^{-1}\sigma$
induces an isometry of $L$ and $\mu \iota(e_{\alpha})=c_\a\iota(e_{\mu \alpha})$
for $\alpha\in L$ and some constant $c_\a\in \C.$ 
Let $\tau\in O(\hat L)$ such that $\bar\tau=\mu.$   
Then $\mu\tau^{-1}$ is the identity when restricted to $\mathfrak{h}(-1).$ 
By Lemma \ref{lemma:2.6} there exists $\psi\in N$ such that
$\mu\tau^{-1}=\psi.$ That is, $\sigma=\phi\psi\tau\in N\cdot O(\widehat{L})$
as $\phi\psi\in N.$ This shows that $\Aut(V_L)=N\cdot O(\hat L).$

In order to see $\Hom (L,\Z/2\Z)$ is a subgroup of $N$ we recall from
[FLM] on how $\Hom (L,\Z/2\Z)$ is embedded into $\Aut(\hat L).$ 
Let $\lambda\in \Hom (L,\Z/2\Z).$ The corresponding automorphism
of $\hat L$ is given by: $a\mapsto a(-1)^{\lambda(\bar a)}$ for
$a\in\hat L.$ As in the proof of Lemma \ref{lemma:2.7} there 
exists $h\in \mathfrak{h}$ such that $e^{\<h,\alpha\>}=(-1)^{\lambda(\a)}$
for $\a\in L.$ It is easy to see that the $\lambda$ and $e^{h(0)}$
give the same automorphism of $V_L$ and thus  $\lambda\in N.$ 

An isomorphism from $\Aut(V_L)/N$ to a quotient group of $O(L)$ follows
from the facts that $\Aut(V_L)/N$ is isomorphic to 
$O(\hat L)/N\cap O(\hat L),$ that
$N\cap O(\hat L)$ contains $\Hom (L,\Z/2\Z)$ and that
$O(\hat L)/\Hom (L,\Z/2\Z)$ is isomorphic to $O(L).$

\section{$\theta$-twisted modules}
\setcounter{equation}{0}

Let $V_L$ be as before. Fix an involution $\theta$ of $\hat L$ of
order 2 such that $\bar \theta=-1$ on $L.$ For example, one can
define $\theta(a)=a^{-1}(-1)^{\<\bar a,\bar a\>/2}$ for $a\in \hat L.$
Then $\theta$ is an automorphism of $V_L$ and $M(1).$ We determine
twisted Zhu's algebra $A_\theta(V_L)$ and $A_\theta(M(1))$ explicitly
in this section. We also classify the irreducible $\theta$-twisted
modules for these vertex operator algebras by using $A_\theta(V_L)$ 
and $A_\theta(M(1))$ and prove the $\theta$-rationality of $V_L.$ 
 
Beginning with a arbitrary vertex operator algebra $V$ and an
automorphism $g$ of finite order $g$ in Subsection 3.1, we recall
the definitions of admissible $g$-twisted modules for $V$ and
twisted Zhu's algebra $A_g(V)$ and related results from [DLM3]. The
$g$-rationality of $V$ is established under certain assumption
on the structure of $A_g(V).$ In Subsection 3.2 we give the construction
of $\theta$-twisted modules for $V_L$ following [FLM]. In Subsection
3.3 we obtain a finite spanning set of $A_\theta(V_L)$ indexed by
$L/2L$ and prove that $A_\theta(M(1))$ is isomorphic to the
one-dimensional associative algebra $\C.$ 
We show in Subsection 3.4 that $A_\theta(V_L)$ is isomorphic to
a quotient of the group algebra $\C[\hat L/K]$ where $K$ is a 
subgroup of $\hat L$ (defined in Subsection 3.2) isomorphic to
$2L.$ Using this result we prove that every irreducible $\theta$-twisted
$V_L$-module is isomorphic to $V_L^{T_\chi}$ for some $\hat L/K$-modules
$T_\chi$ (see Subsection 3.3 for the definition of $T_\chi$). 
The $\theta$-rationality of $V_L$ is an easy corollary of
a general result on $g$-rationality in Subsection 3.1. 

\subsection{Twisted modules and twisted Zhu algebras}

Let $V$ be a vertex operator algebra  and
$g$ an automorphism of $V$ of finite order $T.$  
Decompose $V$ into $g$-eigenspaces: $V=\bigoplus_{r\in \Z/T\Z}V^r$
where $V^r=\{v\in V|gv=e^{-2\pi ir/T}v\}$.

An admissible $g$-twisted $V$-module (cf. [DLM3], [Z]) is a 
$\frac{1}{T}\Z$-graded vector space
$$M=\sum_{n=0}^{\infty}M(\frac{n}{T}),$$
with top level 
$M(0)\ne 0,$ equipped 
with a linear map
$$\begin{array}{l}
V\to (\End\,M)\{z\}\\
v\mapsto Y_M(v,z)=\sum_{n\in\Q}v_nz^{-n-1}\ \ \ (v_n\in
\End\,M)
\end{array}$$
satisfying the following conditions;
for all $0\leq r\leq T-1,$ $u\in V^r$, $v\in V,$ 
$w\in M$,
\begin{eqnarray*}
& &Y_M(u,z)=\sum_{n\in \frac{r}{T}+\Z}u_nz^{-n-1}, \label{1/2}\\ 
& &u_nw=0\ \ \                                  
\mbox{for}\ \ \ n\gg 0,\label{vlw0}\\
& &Y_M({\bf 1},z)=1,\label{vacuum}
\end{eqnarray*}
\[
\begin{array}{c}
\displaystyle{z^{-1}_0\delta\left(\frac{z_1-z_2}{z_0}\right)
Y_M(u,z_1)Y_M(v,z_2)-z^{-1}_0\delta\left(\frac{z_2-z_1}{-z_0}\right)
Y_M(v,z_2)Y_M(u,z_1)}\\
\displaystyle{=z_2^{-1}\left(\frac{z_1-z_0}{z_2}\right)^{-r/T}
\delta\left(\frac{z_1-z_0}{z_2}\right)
Y_M(Y(u,z_0)v,z_2)},
\end{array}
\]
where $\delta(z)=\sum_{n\in\Z}z^n$ and
all binomial expressions are to be expanded in nonnegative
integral powers of the second variable;
$$u_mM(n)\subset M(\wt(u)-m-1+n)$$
if $u$ is homogeneous. 

A $g$-{\em twisted $V$-module} is
an admissible $g$-twisted $V$-module $M$ which carries a 
$\C$-grading  by weight. That is, we have
$$M=\coprod_{\lambda \in{\C}}M_{\lambda} $$
where $M_{\lambda}=\{w\in M|L(0)w=\lambda w\}.$ Moreover we require that 
$\dim M_{\lambda}$ is finite and for fixed $\lambda,$ $M_{\lambda+\frac{n}{T}}=0$
for all small enough integers $n.$ 

$V$ is called $g$-{\em rational} if any admissible $g$-twisted module
for $V$ is a direct sum of irreducible $g$-twisted $V$-modules.

Nest we review twisted Zhu's algebra [DLM3]. For homogeneous $u\in V^r$ and $v\in V$ we  define
\begin{eqnarray}\label{g2.2}
u\circ_g v=\Res_{z}
\frac{(1+z)^{\wt(u)-1+\delta_{r}+\frac{r}{T}}}{z^{1+\delta_r}}Y(u,z)v
\end{eqnarray}
where $\delta_{r}=1$ if $r=0$ and  $\delta_{r}=0$ if $r\ne 0$ and
where, here and below, $(1+z)^{\alpha}$ for $\alpha\in\C$ is to
be expanded in nonnegative integer powers of $z.$ Let $O_g(V)$ be the linear span of all $u\circ_g v$ and define the linear space
$A_g(V)$ to be the quotient $V/O_{g}(V).$

\begin{lemma}\label{lemma:3.2}
(1) Assume that $u\in V^{r}$ is homogeneous,
$v\in V$ and $m\ge n\ge 0.$ Then 
$$\Res_{z}\frac{(1+z)^{{\wt}(u)-1+\delta_{r}+\frac{r}{T}+n}}{z^{m+\delta_{r}+1}}Y(u,z)v\in O_{g}(V).$$

(2) $V^r$ is a subspace of $O_g(V)$ if $r\ne 0.$ 
\end{lemma}

Let us define the second product $*_g$ on $V$.
For $r,u$ and $v$ as
above, set
\begin{equation}\label{a5.1}
u*_gv=\left\{
\begin{array}{ll}
\Res_zY(u,z)\frac{(1+z)^{\wt(u)}}{z}v
 & {\rm if}\ r=0\\
0  & {\rm if}\ r>0.
\end{array}\right.
\end{equation} 
Extend this to a bilinear product  on $V.$

Let $M$ is an admissible $g$-twisted $V$-module. Define the ``vacuum
space'' $\Omega(M)$ of $M$ to be 
$$\{w\in M|u_nw=0, n>\wt(u)-1, n\in \frac{1}{T}\Z\}.$$
Then $M(0)\subset \Omega(M)$ and the weight zero operator $o(v),$ 
which is equal to $v_{\wt (v)-1}$ when $v$ is homogeneous, 
acts on $\Omega(M).$  

\begin{theorem}\label{ta}{\rm [DLM3]} (1) The product $*_g$ induces the structure of an associative algebra on
$A_g(V)$.

(2) For any admissible $g$-twisted $V$-module $M,$
$v+O_g(V) \mapsto \End \Omega(M)$ gives an $A_g(V)$-module
structure on $\Omega(M)$ and $M(0)$ is a submodule of $\Omega(M).$

(3) $M\mapsto M(0)$ gives a bijection between the equivalence
classes of irreducible $g$-twisted admissible $V$-modules and
the equivalence
classes of simple $Ag(V)$-modules.

(4) If $V$ is $g$-rational,  $A_g(V)$ is a finite-dimensional
semisimple algebra.

\end{theorem}

The following theorem is new.

\begin{theorem}\label{tnew} Let $V$ be a vertex operator algebra and $g$ 
an automorphism of order $T$ such that $A_g(V)$ is a finite-dimensional
semisimple algebra and $\o+O_g(V)$ acts on the all irreducible
modules as the same constant $\lambda.$ Then $V$ is $g$-rational.
\end{theorem}

\begin{proof} By definition of $g$-rationality we need to prove that
any admissible $g$-twisted $V$-module $M$ is completely reducible.

First we show that if $M(0)$ is an irreducible $A_g(V)$-module and
$M$ is generated by $M(0)$ then $M$ is irreducible. In this case
$M=\oplus_{n\in \frac{1}{T}\Z,n\geq 0}M_{\lambda+n}.$ If $M$ is 
not irreducible it contains a submodule $W=\oplus_{n\geq n_0>0}W_{\lambda+n}$
for some $n_0>0$ as $M(0)$ is an irreducible $A_g(V)$-module 
where $W_{\mu}=M_{\mu}\cap W$ and $W_{\lambda+n_0}\ne 0.$ Then
$W_{\lambda+n_0}$ is a module for $A_g(V)$ on which $\o+O_g(V)$ acts
as $\lambda+n_0\ne \lambda.$ Since $A_g(V)$ is
a finite-dimensional semisimple algebra, there exists an irreducible
submodule of $W_{\lambda+n_0}$ on which $\omega+O_g(V)$ acts
as $\lambda+n_0\ne \lambda.$  This is a contradiction.

Now we consider general $M.$ Let $X$ be the submodule of $M$ generated by
$\Omega(M).$ Since $\Omega(M)$ is a direct sum of irreducible $A_g(V)$-module,
$X$ is completely reducible by the argument above. Note the $X$ is the
maximal completely reducible submodule of $M.$  We claim that $X=M.$

If $\bar M=M/X\ne 0.$ Then $\Omega(\bar M)\ne 0.$ We can assume that $\bar M$
is completely reducible (otherwise replace $\bar M$ by the submodule
generated by $\Omega(\bar M)$). Thus any vector of $M$ is contained
in a finite-dimensional $L(0)$-stable subspace and $M$ is a direct
sum of generalized eigenspaces for $L(0):$
$$M=\oplus_{n\in \frac{1}{T}\Z,n\geq 0}M_{\lambda+n}$$
where $M_{\mu}$ is the generalized eigenspace with the eigenvalue 
$\mu.$ Since $M_{\lambda}$ is a completely reducible $A_g(V)$-module,
the submodule $W$ generated by $M_{\lambda}$ is completely reducible. If $M\ne
W$ then $\Omega(M/W)$ contains a $A_g(V)$-submodule on which $\o+O_g(V)$ 
acts as a constant which is greater than $\lambda.$  Thus
$M$ is completely reducible and $M=X.$ This is a contradiction.
\end{proof}

\subsection{The $\theta$-twisted modules}

We return to $V_L.$ Recall that $\theta$ is an automorphism of $\hat L$
and $V_L$ such that 
$\theta(a)=a^{-1}(-1)^{\<\bar a,\bar a\>/2}$ for $a\in \hat L.$
Define $K = \{\theta(a)a^{-1}\,|\,a\in \widehat{L}\}.$ Then
$\bar{K} = 2L$. Also set $R = \{\alpha\in L\,|\,
\langle\alpha,L\rangle\subset 2\Z\}\supset 2L$.
Then the pull back $\widehat{R}$ of $R$ in $\widehat{L}$ is in the center of 
$\widehat{L}$ and $K$ is a subgroup of $\widehat{R}$.

The following proposition can be found in [FLM] (Proposition 7.4.8).
\begin{proposition}\label{p3.1}
There are exactly $|R/2L|$ central characters
$\chi:\widehat{R}/K\longrightarrow \C^\times$ of $\widehat{L}/K$ such that
$\chi((-1)K) = -1$. For each such $\chi$, there is a unique (up to equivalence)
irreducible $\widehat{L}/K$-module $T_\chi$ with central character $\chi$
and every irreducible $\widehat{L}/K$-module on which $(-1)K$ acts
$-1$ is equivalent to ones of these. In particular, viewing $T_\chi$ as an
$\widehat{L}$-module, $\theta(a) = a$ on $T_\chi$ for $a\in \widehat{L}$.
\end{proposition}

For each $T_\chi$, we define a twisted space
\[
V_L^{T_\chi} = M_{\Z+\frac{1}{2}}(1)\otimes T_\chi\cong
S(\widehat{\mathfrak{h}}_{\Z+\frac{1}{2}}^-)\otimes T_\chi
\]
where $M_{\Z+\frac{1}{2}}(1)$ is the module for the Heisenberg algebra
$\widehat{\mathfrak{h}}_{\Z+\frac{1}{2}}=\mathfrak{h}\otimes 
t^{1/2}\C[t,t^{-1}]\oplus \C c$ defined by
\[
M_{\Z+\frac{1}{2}}(1)=
U(\widehat{\mathfrak{h}}_{\Z+\frac{1}{2}})
\otimes_{U(\widehat{\mathfrak{h}}_{\Z+\frac{1}{2}}^+\oplus \C\,c)}\C.
\]
Here $\widehat{\mathfrak{h}}_{\Z+\frac{1}{2}}^+=
\mathfrak{h}\otimes t^{1/2}\C[t]$ acts trivially on $\C$ and $c$ acts
as 1.

The $\theta$-twisted vertex operators 
$Y_\theta(\iota(a),z)$ for $a\in \hat L$ are defined to be
\[
Y_\theta(\iota(a),z)
= 2^{-\langle\bar a,\bar a\rangle}
\NO e^{\int \bar a(z)}\NO a
z^{-\langle\bar a,\bar a\rangle/2}
\]
where $\alpha(z) = \sum_{n\in\Z+\frac{1}{2}}\alpha(n)z^{-n-1}$ and
$\a(n)=\a\otimes t^n$ (see [FLM] for the detail).
For $v=\alpha_1(-n_1)\cdots \alpha_k(-n_k)\otimes \iota(a) \in V_L$ homogeneous, set
\[
W_\theta(v,z)
= \NO
\partial^{(n_1-1)}\alpha_1(z)\cdots\partial^{(n_k-1)}\alpha_k(z)
Y_\theta(\iota(a),z)
\NO
\]
where $\partial^{(n)}=\frac{1}{n!}\left(\frac{d}{dz}\right)^n.$ 
Now we extend this to all $v\in V_L$ by linearity.
Recall that $\{h_1,\dots,h_d\}$ be an orthonormal basis of $\mathfrak{h}$
and define constants $c_{mn}\in\Z$ for $m,n\geq 0$ by
the formula
\[
\sum_{m,n\geq 0}c_{mn}x^my^n = -\log{
\left(
\frac{(1+x)^{\frac{1}{2}}+ (1+y)^{\frac{1}{2}}}{2}
\right)}.
\]
Set
\[
\Delta_z = \sum_{m,n\geq 0}\sum_{i=1}^d
c_{mn}h_i(m)h_i(n)z^{-m-n}.
\]
We finally define twisted vertex operator
$Y_{\theta}(v,z)$ for $v\in V_L$ as
\[
Y_\theta(v,z)= W_\theta(e^{\Delta_z}v,z).
\]

The following theorem was proved in [FLM] (see also [DL2]).
\begin{theorem}\label{t3.2} (1) For each $T_\chi$, the space $(Y_\theta,V_L^{T_\chi})$
is an irreducible $\theta$-twisted module for the vertex operator algebra
$V_L$.

(2) $M_{\Z+\frac{1}{2}}(1)$ is an irreducible $\theta$-twisted
$M(1)$-module.
\end{theorem}

Let $W$ be an $\theta$-stable subspace of $V_L.$ Decompose 
$W$ into eigenspaces with respect to
the action of $\theta:$
\[
W = W^+\oplus W^-
\]
where
\[
W^{\pm} = \{v\in W\,|\, \theta(v) = \pm v\}.
\] 

\begin{lemma}\label{lemma:3.3} Let $W=M(1)$ or $V_L.$ 
Suppose $u\in W^-$ and $\wt(u) = 1$. Then
\[
\sum_{j=0}^\infty\binom{\frac{1}{2}}{j}u_{j-m-1}v\in O_\theta(W).
\]
\end{lemma}
\begin{proof}
Take $n=0$ in Lemma \ref{lemma:3.2} (1), then we see
\[
\Res_z\frac{(1+z)^{\frac{1}{2}}}{z^{1+m}}Y(u,z)v\in O_\theta(W).
\]
Substituting $ (1+z)^{\frac{1}{2}} = \sum_{j=0}^\infty\binom{\frac{1}{2}}{j}
z^j$ into the above, we get the result.
\end{proof}

\subsection{A spanning set of $A_\theta(V_L)$}

\begin{lemma}\label{lemma:3.4} Let $W=M(1)$ or $V_L.$ 
Let $\alpha\in \mathfrak{h}$ and $m\in\Z_{\geq 0}$. Then for any $v\in W$,
there exist scalars $c_{m,j}\in\C$ for $j=0,1,2,\dots$ such that
\[
\alpha(-m-1)v \equiv \sum_{j=0}^\infty c_{m,\,j}\alpha(j)v
\quad \mod O_\theta(W).
\]
\end{lemma}
\begin{proof} We prove the lemma by induction on $m.$ 
Since $\alpha\in W^-$, setting $m=0$ in Lemma \ref{lemma:3.3} gives
\[
\sum_{j=0}^\infty\binom{\frac{1}{2}}{j}\alpha(j-1)v
\in O_\theta(W).
\]
Namely, 
\[
\alpha(-1)v\equiv -\sum_{j=0}^\infty\binom{\frac{1}{2}}{j+1}
\alpha(j)v\quad \mod O_\theta(W).
\]
Then $m=0$ case is done. 

Let $m>0$ and assume that the lemma is true
for all nonnegative integers strictly less than $m$. 
By Lemma \ref{lemma:3.3},
we see
\[
\alpha(-m-1)v \equiv
-\sum_{j=0}^\infty\binom{\frac{1}{2}}{j+1}\alpha(j-m)v
\quad
\mod O_\theta(W).
\]
Using the induction hypothesis completes the proof of the lemma.
\end{proof}

Consider an element
\[
\alpha_1(-n_1)\cdots \alpha_k(-n_k)\otimes \iota (a)
\]
for $\alpha_i\in L,$ $n_1\geq \cdots \geq n_k>0$ and
$a\in\hat L.$ We say the vector $v$ has the length $k$ and
denote it by $\ell(v)$.

\begin{lemma}\label{lemma:3.5}
Let $v\in M(1)\otimes \iota(a)$ for $a\in \hat L.$ Then 
$v \equiv c_v \iota(a) \mod O_\theta(V_L)$
for some scalar $c_v\in\C$. Moreover, if $a=1$ then 
$v \equiv c_v  \mod O_\theta(M(1)).$
\end{lemma}

\begin{proof}
We prove the lemma by induction on the length of $v$.
Clearly it is enough to show the lemma for homogeneous elements.
If $\ell(v)= 0$, the lemma is obviously true.
Suppose that the lemma is true for any homogeneous elements
whose lengths are strictly less than $\ell(v)=k>0$.
Then $v$ has the form
$v = \alpha(-m-1)v'$
for some $m\in\Z_{\geq 0}$ and $\alpha\in L$ where 
$v'\in M(1)\otimes \iota(a)$ and $\ell(v') = k-1$. 
Then by Lemma \ref{lemma:3.4}, we see
\[
v \equiv \sum_{j=0}^\infty c_{m,j}\alpha(j)v'
\quad\mod O_\theta(V_L).
\]
Since $\ell(\alpha(j)v')<t$ for all $j\geq 0$, by induction hypothesis, 
we prove the lemma.

If $a=1$ the result is clear now. 
\end{proof}

From Theorems \ref{ta}, \ref{t3.2} and Lemma \ref{lemma:3.5} we immediately have
\begin{corollary} $A_\theta(M(1))$ is isomorphic to $\C$ and
$M_{\Z+\frac{1}{2}}(1)$ is the unique irreducible $\theta$-twisted
$M(1)$-module. 
\end{corollary}

Recall that $e:L\to \hat L$ is a section. We may and do  assume that
$e_{2\alpha}\in K.$ 

\begin{lemma}\label{lemma:3.6} Let $\alpha\in L.$ 

(1) We have 
$\iota(e_{2\alpha}) \equiv 2^{-4\langle\alpha,\alpha\rangle}\1$ modulo
$O_\theta(V_L).$ 

(2) For any $v\in M(1)\otimes \iota(e_{-\a}),$ 
$v \equiv c_v \iota(e_\a) \mod O_\theta(V_L)$ for some scalar $c_v\in\C$.

\end{lemma}
\begin{proof} (1)
Let $E^\alpha = \iota(e_\alpha)+\iota(\theta e_{\alpha})$. Then clearly $E^\alpha\in V_L^+$
and $\wt(E^\alpha) = \langle\alpha,\alpha\rangle/2$.
Then by Lemma \ref{lemma:3.2} (1), we see that for $m\in\Z_{\geq 0}$ and
$v\in V_L$
\[
E^\alpha_{-m-2}v\equiv
-\sum_{j=1}^{\frac{1}{2}\<\a,\a\>}\binom{\frac{1}{2}\<\a,\a\>}{j}
E^\alpha_{j-m-2}v\quad
\mod O_\theta(V_L).
\]
Let $m = \langle\alpha,\alpha\rangle -1$ and $v = \iota(e_\alpha)$.
Note that
$\iota(e_\alpha)_{j-\langle\alpha,\alpha\rangle-1}\iota(e_\alpha) 
=\delta_{j0}\epsilon(\alpha,\alpha)\iota(e_{2\alpha})$ for
$j\geq 0.$ 
Then we see
\begin{eqnarray*}
& &\epsilon(\alpha,\alpha)\iota(e_{2\alpha})
+\iota(\theta e_{\alpha})_{-\langle\alpha,\alpha\rangle-1}\iota(e_\alpha)\\
& &\ \ \ \ \ \equiv -\sum_{j=1}^{\frac{1}{2}\<\a,\a\>}
\binom{\frac{1}{2}\<\a,\a\>}{j}
\iota(\theta e_{\alpha})_{j-\langle\alpha,\alpha\rangle-1}\iota(e_\alpha)
\quad\mod O_\theta(V_L).
\end{eqnarray*}
Since $\iota(\theta e_{\alpha})_{j}\iota(e_\alpha)\in M(1)$ for all $j$ we see from Lemma \ref{lemma:3.5} that
$\iota(e_{2\alpha}) \equiv c_\alpha\1$
for some scalar $c_\alpha\in\C.$ 
 Finally evaluation of this relation on
top levels of all modules $V_L^{T_\chi}$ gives 
$c_\alpha = 2^{-4\langle\alpha,\alpha\rangle}$. 

For (2) we first use Lemma \ref{lemma:3.5} to get a constant $d_\alpha$ such
that $v \equiv d_v \iota(e_{-\a})$ modulo  $O_\theta(V_L).$ Note that 
$\iota(e_{-\a})=\pm \iota(\theta e_\a).$ Since $\iota(e_\a)-\iota(\theta e_\a)\in V_L^-\subset O_\theta(V_L)$ we conclude $\iota(e_{-\a})\equiv \pm
\iota(e_\a)$ modulo $O_\theta(V_L),$ which finishes the proof of (2). 
\end{proof}

Let $L=\cup_{i\in L/2L}(2L+\beta_i)$ be the coset decomposition of $2L$ in
$L.$ Then 
\begin{proposition}\label{proposition:3.7}
$A_\theta(V_L)$ is spanned by 
\[
\{ \iota(e_{\b_i})+O_\theta(V_L)\,|i\in L/2L\}.
\]
\end{proposition}

\subsection{Algebraic structure of $A_\theta(V_L)$}

Define the elementary Shur polynomials $p_r(x_1,x_2, \dots),\,r\in\Z_{\geq 0}$
in variables $x_1,x_2,\dots$ by the equation
\[
\exp{\left(
\sum_{n=1}^\infty \frac{x_n}{n}y^n
\right)
}
=
\sum_{r=0}^\infty p_r(x_1,x_2,\dots)y^r.
\]

Let $\alpha,\beta\in L$. Then
\begin{align*}
 Y(\iota(e_\a),z)\iota(e_\beta)&=\epsilon(\alpha,\beta)
z^{\langle\alpha,\,\beta\rangle}
\exp{\left(
\sum_{n=1}^\infty\frac{\alpha(-n)}{n}z^n
\right)\iota(e_{\alpha+\beta})
}\\
&=\epsilon(\alpha,\beta)\sum_{r=0}^\infty
p_r(\alpha(-1),\alpha(-2),\dots)\iota(e_{\alpha+\beta})
z^{r+\langle\alpha,\,\beta\rangle}.
\end{align*}
Thus 
\begin{align*}
\iota(e_\a)_{i-1}\iota(e_\beta) & =0\quad (i\geq -\langle\alpha,\beta\rangle+1)\\
\iota(e_\a)_{i-1}\iota(e_\beta) & =\epsilon(\alpha,\beta)
p_{-i-\langle\alpha,\beta\rangle}(\alpha(-1),\alpha(-2),\dots)
\iota(e_{\alpha+\beta})\quad
(i\leq -\langle\alpha,\beta\rangle).
\end{align*}

By using these formulas, we will calculate $\iota(e_\a)*_\theta \iota(e_\b)$
as follows. Set 
\[
E^\a = \frac{1}{2}(\iota(e_\a)+\iota(\theta e_{\alpha})),
\quad F^\alpha = \frac{1}{2}(\iota(e_\a)-\iota(\theta e_{\alpha})).
\]
Then $E^\a\in V_L^+,$ $F^\alpha\in V_L^-$ and
$\iota(e_\a) = E^\alpha + F^\alpha$. From the definition of 
the $*_\theta$ operation and Lemma \ref{lemma:3.2}, we have
\[
V_L^+*_\theta V_L^+\subset V_L^+,\quad
V_L^+*_\theta V_L^-\subset V_L^-\subset O_\theta(V_L),\quad
V_L^-*_\theta V_L = 0.
\]
Thus
\[
\iota(e_\a)*_\theta \iota(e_\b)\equiv
\Res_z\frac{(1+z)^{\frac{1}{2}\<\a,\a\>}}{z}Y(E^\a,z)E^\beta
\quad
\mod O_\theta(V_L).
\]
More exactly, by using component operators,
\[
\iota(e_\a)*_\theta \iota(e_\b)\equiv
\sum_{i=0}^{\frac{1}{2}\<\a,\a\>}\binom{\frac{1}{2}\<\a,\a\>}{i}
E^\a_{i-1}E^\beta\quad \mod O_\theta(V_L)
\]
where
\begin{align*}
E^\alpha_{i-1}E^\beta
&=\frac{1}{4}(\iota(e_\alpha)_{i-1}\iota(e_\b) +\iota(\theta e_{\alpha})_{i-1}\iota(\theta e_{\beta}))\\
&\quad +\frac{1}{4}(\iota(e_\alpha)_{i-1}\iota(\theta e_{\beta})+
\iota(\theta e_{\alpha})_{i-1}\iota(e_{\beta})).
\end{align*}

Let us suppose $\langle\alpha,\beta\rangle<0$. Then
\[
\iota(e_\a)_{i-1}\iota(\theta e_\b)=
\iota(\theta e_\a)_{i-1}\iota(e_\b)= 0
\]
for $i\geq 0$. Therefore
\[
\iota(e_\a)*_\theta \iota(e_\b) 
\equiv \frac{1}{4}\sum_{i=0}^{\frac{1}{2}\<\a,\a\>}
\binom{\frac{1}{2}\<\a,\a\>}{i}
(\iota(e_\a)_{i-1}\iota(e_\b) + \iota(\theta e_\a)_{i-1}
\iota(\theta e_\b))\!\mod O_\theta(V_L).
\]
Using Lemmas \ref{lemma:3.5} and \ref{lemma:3.6} 
produce scalar $c_{\alpha,\,\beta}\in\C$ such that
$\iota(e_\a)*_{\theta}\iota(e_\b) \equiv  
c_{\alpha,\beta}\iota(e_{\alpha+\beta})
\mod O_\theta(V_L)$. Evaluating this on the top levels of modules
$V_L^{T_\chi}$ as in Lemma \ref{lemma:3.6},  we see 
$c_{\alpha,\,\beta} = \epsilon(\alpha,\beta)2^{2\langle\alpha,\,\beta\rangle}$.
In particular, $c_{\alpha,\,\beta}\neq 0$.

If $\langle\alpha,\beta\rangle >0,$ a similar argument shows that
$\iota(e_\a)*_{\theta}\iota(e_\b) \equiv c_{\alpha,\,\beta}\iota(e_{-\alpha+\beta})$
for some scalar $c_{\alpha,\,\beta}\in \C$.
By Lemma \ref{lemma:3.6}, we see
$\iota(e_{-2\alpha})*_\theta \iota(e_{\alpha+\beta}) \equiv 
2^{-4\langle\alpha,\alpha\rangle}\iota(e_{\alpha+\beta})$.
On the other hand, since $\langle -2\alpha,\alpha+\beta\rangle
= -2\langle\alpha,\alpha\rangle-2\langle\alpha,\beta\rangle<0$,
we have 
$\iota(e_{-2\alpha})*_\theta \iota(e_{\alpha+\beta}) \equiv c_{-2\alpha,\alpha+\beta}\iota(e_{-\alpha+\beta})$ for some nonzero constant $c_{-2\alpha,\alpha+\beta}.$
Thus $\iota(e_\a)*_{\theta}\iota(e_\b)\equiv d_{\alpha,\beta}\iota(e_{\alpha+\beta})$
with some scalar $d_{\alpha,\,\beta}\in \C$. Again as before we have
$d_{\alpha,\,\beta}=\e(\a,\b)2^{2\langle\alpha,\,\beta\rangle}.$

For $\alpha\in L$ set $u_\a=2^{\langle\alpha,\a\>}\iota(e_\a).$ Recall
that $R$ is a sublattice of $L.$

\begin{proposition}\label{p3.11} (1)
For any $\alpha,\beta\in L$, 
\[
u_\a*u_\b=\e(\a,\b)u_{\a+\b} \quad \mod O_\theta(V_L).
\]

(2) $u_\a+ O_\theta(V_L)$ is 
in the center of $A_\theta(V_L)$ if and only if $\alpha\in R.$ 
\end{proposition}
\begin{proof} (1) is what we have proved. For (2) we compare 
$$u_\a*u_\b=\e(\a,\b)u_{\a+\b}\quad \mod O_\theta(V_L)$$
and 
$$u_\b*u_\a=\e(\b,\a)u_{\a+\b}\quad \mod O_\theta(V_L).$$
So $u_\a+ O_\theta(V_L)$ is in the
center if and only if $\e(\a,\b)\e(\b,\a)=(-1)^{\<\a,\b\>}=1$
or if and only if $\a\in R.$
\end{proof}

Motivated by Proposition \ref{p3.11} we define a map $f: \C[\hat L/K]
\to A_\theta(V_L)$ such that $f(e_{\alpha}K)=u_\alpha$
and $f((-1)K)=-1.$ Then from Proposition \ref{p3.11}, this is an onto 
algebra homomorphism and $(-1)K+1$ is in the kernel. Let $I$ be
the two sided ideal of  $\C[\hat L/K]$ generated by $(-1)K+1.$ Then
$f$ induces an onto algebra homomorphism $\bar f:  \C[\hat L/K]/I\to 
A_\theta(V_L).$ 

Our second main theorem in this paper is:
\begin{theorem} (1) $A_\theta(V_L)$ is an finite-dimensional
semisimple algebra isomorphic to $\C[\hat L/K]/I.$

(2) The irreducible $\theta$-twisted modules up to isomorphism
are $V_L^{T_\chi}$ for all $T_{\chi}$ given in Proposition \ref{p3.1}.
Moreover, every irreducible admissible $\theta$-twisted $V_L$
is an ordinary $\theta$-twisted $V_L$ module.

(3) $V_L$ is $\theta$-rational.

\end{theorem}

\begin{proof} (1) First note from Proposition \ref{p3.1} and Theorem 
\ref{t3.2} that $\C[\hat L/K]/I$ is a finite-dimensional semisimple
algebra whose 
irreducible modules are precisely those $T_{\chi}$ which are the top levels of
the  irreducible $\theta$-twisted modules $V_L^{T_\chi}$ for 
$V_L.$  So every irreducible module for $\C[\hat L/K]/I$ is an irreducible
module for $A_\theta(V_L)$ and $\bar f$ must be an isomorphism.
 
(2) follows from Theorem \ref{ta} and part (1).

(3) By Theorem \ref{tnew} it is enough to show that $L(0)$ has
the same eigenvalue on all $T_{\chi}.$  In fact it is shown in [FLM]
that $L(0)$ acts on all $T_{\chi}$ as $\frac{1}{16}d$ where $d$ is
the rank of $L.$ 
\end{proof}

\end{document}